\documentclass{llncs}

\usepackage{amssymb}
\usepackage{amsmath}
\usepackage{enumerate}
\usepackage{color}
\usepackage{etoolbox}
\usepackage{graphicx}

\newcommand*{\LargerCdot}{\raisebox{-0.25ex}{\scalebox{1.2}{$\cdot$}}}

\DeclareMathOperator    \aff                    {aff}

\DeclareMathOperator    \argmax         {arg\,max}

\DeclareMathOperator    \conv           {conv}

\DeclareMathOperator    \im                   {im}
\DeclareMathOperator    \intr                   {int}

\DeclareMathOperator    \proj           {proj}
\DeclareMathOperator    \rec                    {rec}

\DeclareMathOperator    \relint         {rel\,int}

\DeclareMathOperator    \verts          {vert}

\AtBeginEnvironment{bmatrix}{\setlength{\arraycolsep}{7pt}}

\newcommand{\R}{\mathbb{R}}

\newcommand{\Z}{\mathbb{Z}}

\newcommand{\ceil}[1]{\lceil #1 \rceil}

\newtheorem{theoremx}[theorem]{Theorem}

\newtheorem{lemmax}[theorem]{Lemma}

\newtheorem{remarkx}[theorem]{Remark}
\newtheorem{observation}[theorem]{Observation}
\spnewtheorem*{proofx}{Proof}{\itshape}{\rmfamily}
\AtEndEnvironment{proofx}{\null\hfill\qed}
\newtheorem{conjecturex}[theorem]{Conjecture}
\newtheorem{examplex}[theorem]{Example}
\newtheorem{definitionx}[theorem]{Definition}

\newcommand{\0}{{\mathbf 0}}
\newcommand{\1}{{\mathbf 1}}

\newcommand{\stack}[2]{\begin{bmatrix} #1 \\ #2 \end{bmatrix}}

\begin{document}

\pagestyle{headings}

\mainmatter 

\title{Mixed Integer Reformulations of Integer Programs and the Affine TU-dimension of a Matrix}

\author{
J\"org Bader\inst{1} 
\and Robert Hildebrand\inst{2} 
\thanks{rhildeb@us.ibm.com}
\and Robert Weismantel\inst{1} 
\and Rico Zenklusen\inst{1} 
}
\institute{Institute for Operations Research, ETH Z\"urich, Switzerland \and   IBM T.J. Watson Research Center}

\maketitle
\begin{center}
\today
\end{center}
 
\begin{abstract}
We study the reformulation of integer linear programs by means of a mixed integer linear program with fewer integer variables.  Such reformulations can be solved efficiently with mixed integer linear programming techniques.
We exhibit examples that demonstrate how integer programs can be reformulated using far fewer integer variables.  
To this end, we introduce a generalization of total unimodularity called the \emph{affine TU-dimension} of a matrix and study related theory and algorithms for determining the affine TU-dimension of a matrix.  We also present bounds on the number of integer variables needed to represent certain integer hulls.


\end{abstract}



\section{Introduction}
\label{sec:Introdcution}

Reformulations of integer programs with linear constraints are common in the integer programming literature.  The main motivation behind these reformulations is that linear programs can be rapidly solved in both theory and practice. 
Often integer programs are reformulated such that the linear programming relaxation is improved, hence creating better bounds for branch-and-bound based algorithms~\cite{50-years-vanderbeck-wolsey}.  These reformulations include Lagrangian relaxation, Dantzig-Wolfe reformulation, and cutting planes.  In a much stronger sense, many combinatorial optimization problems can be formulated exactly using linear inequalities. These formulations typically involve totally unimodular (TU) matrices or totally dual integral systems. 
Since the feasible set of the linear relaxation exactly describes its integer hull in these settings, the optimization problem can be solved by solving the linear relaxation.  See, for instance,~\cite{schrijver_theory_1986}.  

The aim of this work is to instead reformulate integer linear programs by means of a mixed integer linear program with few integer variables. 
We focus on reformulating the feasible region of the integer linear program by mixed integer constraints. For a polyhedron $P = \{ x \in \R^n \mid Ax \leq b\}$, we wish to find a matrix $W \in \Z^{k \times n}$ such that
\begin{equation}
\conv(P \cap \Z^n) = \conv(\{ x \in P \mid Wx \in \Z^k\}).
\label{eq:Wprop}
\end{equation}
Since $W$ is an integral matrix, the inclusion $\subseteq$ in property \eqref{eq:Wprop} is always fulfilled. Furthermore, if $W$ is the $n \times n$ identity matrix, the property \eqref{eq:Wprop} is trivially satisfied.  The objective here is to find an integer matrix $W$ with few rows $k$ that can model the integer hull of the feasible region as in \eqref{eq:Wprop}.    

With such a reformulation at hand, the underlying optimization problem can then be solved using mixed integer linear programming techniques.  From a theoretical point of view, Lenstra~\cite{lenstra_integer_1983} presented an algorithm to solve mixed integer linear programs in polynomial time when the number of integer variables is fixed.  Also in practice, mixed integer linear programs with few integer variables can lead to algorithms with improved running time~\cite{lodi-personal}.

As a centerpiece of this paper, we study the following generalization of total unimodularity that admits reformulations of $P$ for all integral right hand sides $b$.
\begin{definitionx}\label{def:affineTUdecomp}
We say a matrix $A\in \Z^{m\times n}$ admits a \emph{$k$-row affine TU decomposition} if there exist matrices $U \in \Z^{m \times k}$ and $W \in \{0,\pm 1\}^{k \times n}$ such that $\begin{bmatrix} \tilde{A} \\ W \end{bmatrix}$ is totally unimodular, where $\tilde{A}\in\Z^{m\times n}$ is the unique matrix satisfying $A=\tilde{A}+UW$.
The minimum $k$ such that $A$ admits a $k$-row affine TU decomposition is called the \emph{affine TU-dimension of $A$}.  
\end{definitionx}

Observe that $A$ has affine TU-dimension $0$ if and only if $A$ is TU.  Also, the affine TU-dimension of any matrix is at most $n$ since $\tilde A = 0_{m\times n}$, $U = A$ and $W = I_n$ produces a valid affine TU decomposition.  Affine TU decompositions with $k$ rows can admit a model of the feasibility region of related polyhedra with only $k$ integer variables.

\begin{theoremx}\label{thm:TU-decomp-has-int-prop}
 Let $A=\tilde{A}+UW\in\Z^{m\times n}$ with $W\in\{0,\pm 1\}^{k\times n}$ be an affine TU decomposition, $b\in\Z^{m}$ and $\ell \in  (\Z \cup \{-\infty\})^n$, $u \in  (\Z \cup \{\infty\})^n$ with $\ell\le u$.
 
 Then $\conv\left(\left\{ x\in\R^{n}\mid \ell\le x\le u,Ax\le b,Wx\in\Z^{k}\right\} \right)$
is an integral polyhedron. In particular, $P=\{x\in\R^n\mid \ell\le x\le u,Ax\le b\}$ and $W$ satisfy property \eqref{eq:Wprop}.
\end{theoremx}

As an example of the  convenience of Theorem~\ref{thm:TU-decomp-has-int-prop}, the parity polytope has a simple description using one integer variable.

\begin{examplex}[Parity polytope]
The $n$-dimensional (even) parity polytope is the convex hull of all 0-1 vectors in $\R^n$ that have an even cardinality support.  This polytope has exponentially many inequalities~\cite{Jeroslow1975}, but can be described by an extended formulation with only $4n-1$ inequalities~\cite{carrkonjevod2004,Kaibel2013}.   Alternatively, it can be described using one integrality constraint as 
$
P=\conv\left(\left\{x \in [0,1]^n \;\middle\vert\; \tfrac12 \sum_{i=1}^n x_i \in \Z\right\}\right).
$
$P$ is the projection of the polytope
\begin{align*}
Q=&\conv\left(\left\{ \stack{x}{z} \in [0,1]^n \times \R \;\middle\vert\; \sum_{i=1}^n x_i + 2z  = 0, z \in \Z\right\}\right)\\
=&\conv\left(\left\{ \stack{x}{z} \in [0,1]^n \times \R \;\middle\vert\; A\stack{x}{z} \leq b, W\stack{x}{z} \in \Z\right\}\right)
\end{align*}
with $A = \begin{bmatrix} \phantom{-}\1^T_n & 2 \\ -\1^T_n & -2 \end{bmatrix}$, $b = \begin{bmatrix}0\\ 0\end{bmatrix}$ and $W = [\0^T_n, 1]$.  The matrix $A$ admits an affine TU decomposition as $A = \begin{bmatrix} \phantom{-}\1^T_n & 0 \\ -\1^T_n & 0 \end{bmatrix} + \begin{bmatrix}2\\ -2\end{bmatrix} W $. By Theorem~\ref{thm:TU-decomp-has-int-prop}, $Q$ is an integral polyhedron, thus $P$ is also an integral polyhedron.
\end{examplex}

In Section~\ref{sec:TU-decomp_properties}, we expand upon the theory of affine TU decompositions:  we investigate structure of affine TU decompositions and prove Theorem~\ref{thm:TU-decomp-has-int-prop} and related results.

In Section~\ref{sec:TU-decomp_examples}, we give various examples of how understanding the affine TU-dimension can create a mixed integer model with few integer variables.

In Section~\ref{sec:TU-dim_determining}, we focus on computational issues in connection to affine TU decompositions and study the complexity of determining the affine TU-dimension of a matrix.  In particular, we show that it is \textit{NP}-Hard to decide if the (affine) TU-dimension of $A\in\Z^{m\times n}$ is less than $n$. When $k$ and the number $m$ of rows of $A$ is fixed, we give a polynomial time algorithm to determine if the affine TU-dimension of $A$ is equal to $k$. 

In Section~\ref{sec:reformulations-knapsack}, we study mixed integer reformulations for knapsack polytopes. We prove a general lower bound of $\frac{n}{2}$ integrality constraints necessary to achieve property~\eqref{eq:Wprop} for the linear relaxation of knapsack polytopes.
We then give a nonconstructive proof that in every 0-1 knapsack polytope we can replace the integrality on all the $n$ variables by at most $n-2$ integrality constraints to achieve property \eqref{eq:Wprop}. Apart from the added integrality constraints, we use only the original knapsack inequality and the 0-1 bounds on the variables, but no additional inequalities.
In a final example we show the potential power of the addition of linear inequalities in the mixed integer model. We present a class of knapsack polytopes having an exponential sized polyhedral description. For these polytopes, by introducing linearly many additional linear inequality constraints one can replace the integrality constraints on all the variables by a single joint integrality constraint.

\subsection*{Related Work}
The matrices with affine TU-dimension 1 have been called \emph{nearly totally unimodular matrices} in~\cite{Gijswijt2005}. One of several applications of matrices with affine TU-dimension 1 is edge coloring in \textit{nearly bipartite graphs}~\cite{Gijswijt2005}. An undirected graph $G$ is called nearly bipartite if it is not bipartite but one can obtain a bipartite graph by deleting one vertex of $G$~\cite{Eggan198671}. The incidence matrix of a nearly bipartite graph has affine TU-dimension 1. By Theorem~\ref{thm:TU-decomp-has-int-prop}, integer hulls described by nearly totally unimodular matrices can be captured using just one integrality constraint, and hence integer linear optimization on these problems can be done in polynomial time using a mixed integer linear program.  This, however, is not mentioned in~\cite{Gijswijt2005}.

Mixed integer reformulations were studied by Martin~\cite{martin1987} under the name \emph{variable redefinition}.  Martin showed how to reformulate problems that can be solved using dynamic programs.  The main motivation here was to create tighter linear programming relaxations that improve bounds in a branch and bound algorithm.  

Many polynomial time algorithms in combinatorial optimization rely on the method of guessing the value of certain problem-related variables. Depending on the situation, guessing can for instance be done by polynomial enumeration of all possible combinations of values. Often one can interpret this approach as a polyhedral problem with an underlying mixed integer property like in \eqref{eq:Wprop}. As an example, Hassin et al.~\cite{Hassin2004} gave an efficient polynomial time approximation scheme for the constrained minimum spanning tree problem. Their algorithm relies on a partition of the edge set into (logarithmically many) buckets $T_i$. Then it guesses the number of chosen edges $\sum_{e\in T_i}x_e$ in each bucket $T_i$, and constructs an approximate solution from this information. The guessing can be interpreted as solution of the linear relaxation of the polyhedral problem with additional integral constraints $\sum_{e\in T_i}x_e\in\Z$ for all $i$.

Another example of guessing integral values of certain variables is presented in Oriolo et al.~\cite{Oriolo2013}. They provide an exact efficient algorithm for a special case of a network design problem on rings which relies on guessing the right value for an integral variable $b$. Since their integral variable could have an exponential range, a complete enumeration is not efficient. They thus provide another argument, leading to an efficient procedure. The introduction of the integrality constraint $b\in\Z$ gives an alternative conceptually simple way to obtain a polynomial time algorithm.

We would like to mention the close connection of mixed integer reformulations to extended formulations.
Consider the special case that the number of possible values $d=Wx \in \Z^k$ in \eqref{eq:Wprop} is polynomial in $n$. Then one can rewrite 
\[
\conv(\{x\in P\mid Wx\in\Z^k\})=\conv\left(\bigcup_{d\in\Z^k}\{x\in P\mid Wx=d\}\right)
\]
as a linear program of polynomial size in an extended space with a method of Balas~\cite{Balas1998}.
In this special case, a formulation as in~\eqref{eq:Wprop} can serve as a compact certificate for the existence of such an extended formulation.

Lastly, we note that the set $\left\{ x\in\R^{n}\mid Ax\le b,Wx\in\Z^{k}\right\}$ is a projection from a $n+k$ dimensional space. When $W$ is unimodular, as a convenient feature, we can find a change of variables that allows the integrality constraints $Wx\in\Z^k$ to be modeled easily using integral variables without increasing the total number of variables.  To see this, consider the Hermite normal form transformation of a unimodular matrix $W\in\Z^{k\times n}$ with full row rank $k$, that is $W\cdot L=[I_k\; 0_{k\times (n-k)}]$ with $L\in\Z^{n\times n}$ unimodular. Then by the bijective linear transformation $x=Ly$, we have
$$
\max \{c^Tx \mid Ax\le b,Wx\in\Z^{k}\} = \max \{(c^T L) y \mid ALy\le b,\; y_1,\ldots,y_k\in\Z\}.
$$
Using this transformation prevents the need to write the problem in a higher dimension to model the integrality constraints.



\section{Properties of TU decompositions}
\label{sec:TU-decomp_properties}

We will begin with a homogeneous version of affine TU decompositions that is easier to study.  
\begin{definitionx}
We say a matrix $A\in \Z^{m\times n}$ admits a \emph{$k$-row TU decomposition} if there exist matrices $U \in \Z^{m \times k}$ and $W \in \{0,\pm 1\}^{k \times n}$ with $A=UW$ such that $W$ is totally unimodular.
The minimum $k$ such that $A$ admits a $k$-row TU decomposition is called the \emph{TU-dimension of $A$}.  
\end{definitionx}
 An obvious lower bound on the TU-dimension of a matrix is given by its rank, therefore even a TU matrix can have a large TU-dimension. The affine version of TU decompositions can rule out this artifact.


Following 
Theorem~\ref{thm:TU-decomp-has-int-prop}, for a given matrix $A$,  we are interested in finding a $k$-row affine TU decomposition $A=\tilde{A}+UW$ with small $k$ (the number of rows of $W$).  With this purpose in mind, $U$ can be restricted to not contain an all-zero column. Otherwise, we could delete this $r$-th column of $U$ together with the $r$-th row of $W$ and still obtain a valid affine TU decomposition.

Furthermore, we may assume without loss of generality that $W\in\{0,\pm 1\}^{k\times n}$ has full rank $k$. This will allow us to simplify some proofs in what follows, but as the following remark shows it is not a real limitation.

\begin{remarkx}\label{rem:W-full-rank}
Let $A \in \Z^{m \times n}$ and $U \in \Z^{m \times k}$ such that $A = \tilde A + U W$ where $\begin{bmatrix} \tilde{A} \\ W \end{bmatrix}$ is TU and $W\in\{0,\pm 1\}^{k\times n}$ has rank $k'<k$.

For any $I \subseteq \{1, \dots, k\}$, let $W_{I\LargerCdot}$ be the matrix consisting of the rows of $W$ with indices from $I$, and $U_{\LargerCdot I}$ be the matrix consisting of the columns of $U$ with indices from $I$. 

Now fix some $I \subseteq \{1, \dots, k\}$, with $|I| = k'$ such that the TU matrix $W_{I\LargerCdot}$ has rank $k'$.  For any $j \in I^c := \{1, \dots, k\} \setminus I$, consider $W_{j\LargerCdot}$, the $j$-th row of $W$.  Then $(W_{j\LargerCdot})^T=(W_{I\LargerCdot})^Tr^j\in\Z^n$ for some $r^j\in\R^{k'}$. Since $W$ is TU, there exists such a solution $r^j$ with $r^j\in\Z^{k'}$.  

Then $W_{I^c\LargerCdot} = R W_{I\LargerCdot}$ where $R\in \Z^{(k - k') \times k'}$ is the matrix composed of rows $(r^j)^T$ for $j \in I^c$.  
Since $\begin{bmatrix}\tilde A \\ W_{I\LargerCdot}\end{bmatrix}$ is TU, $A = \tilde A + U W = \tilde A + U_{\LargerCdot I} W_{I\LargerCdot} + U_{\LargerCdot I^c} W_{I^c\LargerCdot} = \tilde A + (U_{\LargerCdot I} +  U_{\LargerCdot I^c} R)W_{I\LargerCdot}$
is a $k'$-row affine TU decomposition of $A$.
\end{remarkx}

We now prove Theorem~\ref{thm:TU-decomp-has-int-prop}.   A polyhedron $P\subseteq \R^n$ is called \emph{integral} if $P = \conv(P \cap \Z^n)$.
\begin{proofx}[of Theorem~\ref{thm:TU-decomp-has-int-prop}]
Let $I^1$ be the $r_1\times n$ submatrix of the identity matrix $I_n$ with rows corresponding to the finite indices $\ell_i\in\Z$, and let $I^2$ be the $r_2\times n$ submatrix of $I_n$ with rows corresponding to the finite indices $u_i\in\Z$.\\
Notice that any affine TU decomposition $A=\tilde{A}+UW$ gives rise to an affine TU decomposition
$\hat{A}:=\begin{bmatrix}A\\-I^1\\I^2\end{bmatrix}=\begin{bmatrix}\tilde A\\-I^1\\I^2\end{bmatrix}+\begin{bmatrix}U\\ 0_{r_1\times k}\\ 0_{r_2\times k}\end{bmatrix}W.$
Moreover, with $\hat b:=\begin{bmatrix}b\\-l_{I^1}\\u_{I^2}\end{bmatrix}$ one has $$\conv\left(\left\{ x\in\R^{n}\mid \ell\le x\le u,Ax\le b,Wx\in\Z^{k}\right\} \right)=\conv\left(\left\{ x\in\R^{n}\mid \hat Ax\le b,Wx\in\Z^{k}\right\} \right).$$
Thus, it suffices to prove that for any affine TU decomposition $A=\tilde{A}+UW$, the set $S = \conv\left(\left\{ x\in\R^{n}\mid Ax\le b,Wx\in\Z^{k}\right\} \right)$ is an integral polyhedron.

Observe first that $S$ is a polyhedron.  To see this, let  $$\hat S=\conv\left(\left\{ (x,d)\in\R^{n}\times\Z^k\mid Ax\le b,Wx-d=\0_k\right\}\right).$$  Then $\hat S$ is a polyhedron since it is the mixed integer hull of a rational polyhedron (see, e.g. \cite[Section 16.7]{schrijver_theory_1986}).  Since projections of polyhedra are also polyhedra, the projection  $S = 
\proj_x(\hat S)$ of $\hat S$ onto the $x$ variables is also a polyhedron.

We now show that $S$ is integral.
Consider a decomposition of $A$ as stated. For every fixed $d\in\Z^{k}$,
$b-Ud$ is an integral vector, which implies that 
$$P_{d}:=\left\{ x\in\R^{n}\mid \tilde{A}x\le b-Ud,Wx=d\right\} 
$$ is an integral polyhedron by total unimodularity of $\begin{bmatrix} \tilde{A} \\ W \end{bmatrix}$~\cite[Theorem 19.1]{schrijver_theory_1986}.  Then 
$$
S = \conv\left(\bigcup_{d\in\Z^{k}}P_{d}\right)
   =  \conv\left(\bigcup_{d\in\Z^{k}} \conv(P_{d} \cap \Z^n) \right)
   = \conv\left(\bigcup_{d\in\Z^{k}} (P_{d} \cap \Z^n) \right)
   =  \conv\left(S \cap \Z^n \right),
$$
and thus, $S$ is an integral polyhedron.  The second equality follows from the fact that $P_d$ is an integral polyhedron and the last inequality follows since $S \cap \Z^n = \cup_{d \in \Z^k} P_d \cap \Z^n$.  This is clear since for any $x \in P \cap \Z^n$, $Wx \in \Z^k$ by integrality of $W$, and hence $x \in P_d$ for some $d \in \Z^k$.  

Therefore, $\conv\left(\left\{ x\in\R^{n}\mid Ax\le b,Wx\in\Z^{k}\right\} \right)$ is an integral polyhedron. From the arguments outlined before, also $\conv\left(\left\{ x\in\R^{n}\mid \ell\le x\le u,Ax\le b,Wx\in\Z^{k}\right\} \right)$ is an integral polyhedron.
\end{proofx}

In a slightly restricted setting also a converse of Theorem~\ref{thm:TU-decomp-has-int-prop} holds.

For this, let us first link affine TU decompositions to affine decompositions containing general unimodular matrices. Recall that an integral (not necessarily square) $r\times s$ matrix is called unimodular if it has full row rank $r$ and each nonsingular $r\times r$ submatrix has determinant $\pm 1$. An integral matrix $M\in\Z^{r\times s}$ of rank $r$ is unimodular if and only if the polyhedron $\{x\in\R^s\mid Mx=v,\,x\ge 0\}$ is integral for each $v\in\Z^m$ (see, e.g. \cite[Theorem 19.2]{schrijver_theory_1986}).

Now, let $\begin{bmatrix} A &I_m\\W &0_{k\times m}\end{bmatrix}\in\Z^{(m+k)\times(n+m)}$ have rank $m+k$.
Notice that
$\begin{bmatrix} A &I_m\\W &0_{k\times m}\end{bmatrix}$ is unimodular if and only if 
$\{(x,y)\in\R^{n+m}\mid Ax+y=b,\,Wx=d,\,x\ge 0,\,y\ge 0\}$ is integral for all $b\in\Z^m$ and all $d\in\Z^k$. The latter is true if and only if
$\{x\in\R^{n}\mid Ax\le b,\,Wx=d,\,x\ge 0\}$ is integral for all $b\in\Z^m$ and all $d\in\Z^k$.

For reference we sum up this insight as a remark.

\begin{remarkx}
\label{rem:unimod}
Let $A=\tilde{A}+UW\in\Z^{m\times n}$ be a decomposition of $A$ with only integral matrices.
\begin{enumerate}[(i)]
\item 
The matrix $\begin{bmatrix} A &I_m\\W &0_{k\times m}\end{bmatrix}$ is unimodular if and only if $\left\{ x\in\R^{n}\mid Ax\le b,Wx=d,\,x\ge 0\right\}$ is an integral polyhedron for all $b\in\Z^m$ and all $d\in\Z^k$. 
\item
Moreover, if $\begin{bmatrix} A &I_m\\W &0_{k\times m}\end{bmatrix}$ is unimodular then 
 $\conv\left(\left\{ x\in\R^{n}\mid Ax\le b,Wx\in\Z^{k},\,x\ge 0\right\} \right)$ is an integral polyhedron for each $b\in\Z^{m}$ . 
\end{enumerate}
\end{remarkx}
An anonymous reviewer pointed out to us that the equivalence in (i) above does not hold anymore if one drops the nonnegativity assumption on $x$.

\begin{theoremx}
Let $A\in\Z^{m\times n}$, and let $W\in\Z^{k\times n}$ have rank $k$ such that the polyhedron $\left\{ x\in\R^{n}\mid Ax\le b,Wx=d,\,x\ge 0\right\}$ is integral for all $b\in\Z^m$ and for all $d\in\Z^k$.
Then there exist matrices $U\in\Z^{m\times k}$ and $W'\in\{0,\pm 1\}^{k\times n}$ such that $A=\tilde{A}+UW'$ is an affine TU decomposition. Moreover, for every $b\in\Z^m$, $$\conv\left(\left\{ x\in\R^{n}\mid Ax\le b,Wx\in\Z^k,\,x\ge 0\right\} \right)=\conv\left(\left\{ x\in\R^{n}\mid Ax\le b,W'x\in\Z^k,\,x\ge 0\right\} \right).$$
\end{theoremx}
\begin{proofx}
Since $W$ has rank $k$, we have $n\ge k$. By Remark~\ref{rem:unimod} we know that $\begin{bmatrix} A &I_m\\W &0_{k\times m}\end{bmatrix}$ is a unimodular matrix. We claim that also $W$ is unimodular. For this, consider any nonsingular $k\times k$ submatrix $W_J$ of $W$, where $J$ with $|J|=k$ is the set of columns chosen. Since $B=\begin{bmatrix} A_J &I_m\\W_J &0_{k\times m}\end{bmatrix}$ has rank $m+k$, we have $\det(W_J)=\pm\det(B)\in\{\pm 1\}$. Thus, $W$ is unimodular.

Let $W=[W_1\;\;W_2]$ where we assume without loss of generality that the last $k$ columns of $W$ are the full rank matrix $W_2$. Then we have $$\begin{bmatrix} A &I_m\\W &0_{k\times m}\end{bmatrix}=\begin{bmatrix} A_1&A_2 &I_m\\W_1&W_2 &0_{k\times m}\end{bmatrix}=\begin{bmatrix} I_m&A_2\\ 0_{k\times m}&W_2\end{bmatrix}\cdot\begin{bmatrix} A_1-A_2W_2^{-1}W_1 & 0_{m\times k} &I_m\\W_2^{-1}W_1 &I_k&0_{k\times m}\end{bmatrix}\,.$$
The matrix $\begin{bmatrix} I_m&A_2\\ 0_{k\times m}&W_2\end{bmatrix}$ is unimodular since it is (up to permutation of the columns) a full row rank submatrix of the unimodular matrix $\begin{bmatrix} A &I_m\\W & 0_{k\times m}\end{bmatrix}$. As a product of unimodular matrices, $\begin{bmatrix} A_1-A_2W_2^{-1}W_1 & 0_{m\times k} &I_m\\W_2^{-1}W_1 &I_k& 0_{k\times m}\end{bmatrix}$ is unimodular as well. Moreover, it is well known that a matrix $C\in\Z^{s\times r}$ is TU if and only if $[C\;\;I_s]$ is unimodular (see, e.g. \cite[Section 19.1]{schrijver_theory_1986}), thus $\begin{bmatrix} A_1-A_2W_2^{-1}W_1 \\W_2^{-1}W_1 \end{bmatrix}$ is TU. Trivially, the matrix $\begin{bmatrix} A_1-A_2W_2^{-1}W_1 & 0_{m\times k}\\W_2^{-1}W_1 &I_k\end{bmatrix}$ is TU as well.

In particular, $A=[A_1-A_2W_2^{-1}W_1\;\; 0_{m\times k}] + A_2 [W_2^{-1}W_1\;\;I_k]$ is an affine TU decomposition.
Moreover, $Wx\in\Z^k$ if and only if $[W_2^{-1}W_1\;\;I_k]x\in\Z^k$.
\end{proofx}


\section{Examples of TU decompositions}
\label{sec:TU-decomp_examples}


Given a matrix, in order to find (affine) TU decompositions of it having few rows, one can exploit certain relations between its columns. We first show that any master knapsack polytope
$$P_{\text{Master}} = \conv(\{ x \in \{0,1\}^n \mid  \sum_{i=1}^n i x_i \leq b\})$$
can be modeled with only $\mathcal{O}(\sqrt{n})$ many integrality constraints.  Although there is a well-known polynomial time algorithm for integer linear optimization over $P_{\text{Master}}$ using dynamic programming (see, e.g. \cite[Section 18.5]{schrijver_theory_1986}), this characterization shows that the structure of the resulting problem is less complicated than that of general integer programs.

\begin{examplex}
\label{rem:master_lower_bound}
The affine TU-dimension of 
$
a^T=[1,2,\ldots,n]
$
is $\Theta(\sqrt{n})$.

It follows from Theorem~\ref{thm:TU-decomp-has-int-prop} that $P_{\text{Master}}$ can be modeled with $\mathcal{O}(\sqrt{n})$ many integrality constraints.

Without loss of generality, suppose that $n = \ell^2$ for some $\ell \in \Z_+$.   This is without loss of generality since the affine TU-dimension of $a^T$ as a function of $n$ is monotonically increasing and since $\Theta(\sqrt{n})$, $\Theta(\lceil\sqrt{n}\rceil)$, and $\Theta(\lfloor \sqrt{n} \rfloor)$ are equivalent.  

We begin by exhibiting a matrix $W$ to provide an upper bound.
By defining
\[
W=
\begin{bmatrix}
I_{\ell} & I_{\ell}& I_{\ell} & \ldots & I_{\ell}\\
\0_\ell^T & \mathbf{1}_\ell^T & \0_\ell^T &\ldots & \0_\ell^T\\
\0_\ell^T & \0_\ell^T & \mathbf{1}_\ell^T &  & \0_\ell^T\\
& & &\ddots& \\
  &   &  &  & \mathbf{1}_\ell^T
  \end{bmatrix}
\in\left\{ 0,1\right\} ^{(2\ell-1)\times \ell^{2}}
\]
and
\[
u^T=\begin{bmatrix}
1 & 2 & \ldots & \ell & \ell & 2\ell & \ldots & \ell^2-\ell
\end{bmatrix}
\]
one has a $(2\ell-1)$-row TU decomposition $a^T= u^T W$. Thus $a^T$ has TU-dimension at most $2\ell-1$. Since the first entry of the vector $u$ is 1, we can set $\tilde{a}$ to be the first row of $W$, delete this row in $W$ and the first entry in $u$, and get a $(2\ell-2)$-row affine TU decomposition.

We next explain the lower bound.  
Consider any $k$-row TU decomposition $a^T=u^TW$ with $u\in\Z^k$, and $W\in\{0,\pm 1\}^{k\times n}$ TU. Since the vector $a$ has distinct entries, $W$ needs to have distinct columns. As shown in~\cite[Theorem 4.2]{heller1957}, using slightly different terminology,  a TU matrix with $k$ rows has at most $k^2 + k + 1$ distinct columns (for a more modern approach, see also \cite[Section 21.3]{schrijver_theory_1986}).   Since $W$ has $n$ distinct columns, $k^2 + k + 1 \geq n$.    In particular, $k\ge \sqrt{n}-1$.   A simple argument shows that any $k'$-row affine TU decomposition of $a^T$ gives rise to a $(k'+1)$-row TU decomposition of $a^T$.  Hence, any affine TU decomposition for the master knapsack constraint vector $a^T=[1,2,\ldots,n]$ needs to have at least $\sqrt{n}-2$ rows.
\end{examplex}

\begin{examplex}
The affine TU-dimension of 
$
a^T = [2, 2^2, \dots, 2^n]
$
is $n$.

Consider any TU decomposition of $a^T$ as $a^T = u^TW$ where $W$ has $k$ rows and  full row rank and $u$ has no 0 entries.  Since $W$ is TU, it can be shown that there exist $r^1, \dots, r^{n-k} \in \{0,\pm 1\}^n$ that span the kernel of $W$ (see Lemma~\ref{lem:tu-kernel} in Section~\ref{sec:TU-dim_determining}).  By the decomposition above, $a^T r^i = 0$ for all $i= 1, \dots, n-k$.  But by the structure of $a^T$, observe that $a^T r \neq 0$ for all $r \in \{0,\pm 1\}^n$.  Thus, we must have $k = n$.  
It follows that $W$ is unimodular (and square) and that $W^{-1}\in\mathbb{Z}^{n\times n}$.  Thus we can write $u=(W^{-1})^Ta$. Since the greatest common divisor of the entries in $a^T$ is $2$, $u^T$ has no $\pm 1$ entry.  

Now, consider any affine TU decomposition of $a^T$ as $a^T = \tilde a^T + \bar u^T \bar W$ and suppose the number of rows of $\bar W$ is less than $n$.  Choosing $u^T = [1, \bar u^T]$ and $W = \begin{bmatrix} \tilde a^T \\ \bar W\end{bmatrix}$.  But then $u^T$ and $W$ provide a TU decomposition of $a^T$ with at most $n$ rows in $W$.  This is a contradiction with the above arguments since the first entry of $u^T$ is $1$.  
Thus, the affine TU dimension of $a^T$ is at least $n$.  Since $W$ can be chosen as the identity matrix, the affine TU dimension of $a^T$ is exactly $n$.

\end{examplex}

Next we give an affine TU decomposition for a class of matrices and show how this can be used to model the problem of finding $r$-flows in directed graphs.

\begin{examplex}[Block TU structure]\strut
\begin{enumerate}[(i) ]
\item Let $\left[\begin{matrix}A^i\\ \bar{A}^i
\end{matrix}\right]$ be TU for $i=1,\dots,r$, where $A^i\in\Z^{m_i\times n_i}$, $\bar{A}^i\in\Z^{k_i\times n_i}$ with $M=m_1+\ldots +m_r$.  Let $U^i \in \Z^{k \times k_i}$ for $i=1, \dots, r$ and set
$$
A=
\begin{bmatrix}A^1\\
&A^2\\
&&\ddots\\
&&&A^r\\
U^1 \bar{A^1}& U^2 \bar{A^2}&\ldots&U^r \bar{A^r}
\end{bmatrix}\,.
$$
 Then $A$ has a $(k_1+\ldots +k_r)$-row affine TU decomposition with
$$U=
\begin{bmatrix}
0_{M\times k_1} & 0_{M\times k_2} & \ldots & 0_{M\times k_r}\\
U^1 & U^2 &\ldots&U^r
\end{bmatrix}\text{ and }
W=
\begin{bmatrix}
\bar A^1\\
&\bar A^2\\
&&\ddots\\
&&&\bar A^r
\end{bmatrix}.
$$
\item Consider the task of finding $r$ flows in a directed graph $G=(V,E)$ subject to certain joint capacity constraints as follows. Each flow $f_i\colon E\to\Z$ (for $i=1,\ldots,r$) is allowed to use the arcs $E_i\subseteq E$ to satisfy a demand $d_i$ to be transported from $s_i$ to $t_i$. For each $i$ we demand flow conservation with respect to $f_i$ at all vertices in $V\setminus\{s_i,t_i\}$. Some of the arcs in $E$ are allowed to be shared by several of the flow problems. For each arc $e\in E$, define $I(e)=\{i=1,\ldots,r\mid e\in E_i\}$ to be the set of flow problems involving the arc $e$. We assume that for each arc $e\in E$ we are given a maximum capacity $c_e\in\Z$ that acts as a joint bound $\sum_{i\in I(e)}f_i(e)\le c_e$. 

This task is a variant of the multi-commodity flow problem. In contrast to our variant, in the multi-commodity flow problem one assumes that $I(e)=\{1,\ldots,r\}$ for each $e\in E$.

As is standard for multi-commodity flow problems (see, e.g.~\cite[Chapter 19]{Korte2007}), the above task can be modeled by using integral variables $f_i(e)\in\Z$ for each $e\in E$ and each $i=1,\ldots,r$. Linear constraints $Af\le b$ ensure feasibility of the flows and implement the capacity constraints. Notice that the set of capacity constraints decomposes into the ones for the arcs $e\in E$ with $|I(e)|=1$, and the ones for $e$ with $|I(e)|\ge 2$.

One can rearrange the rows of $A$ in order to get a matrix $A$ as in part~(i). For every $i=1,\ldots,r$, the matrix $A^i$ corresponds to the flow conservation constraints of the flow $f_i$, together with the capacity constraints for $e\in E_i$ with $|I(e)|=1$. Moreover, the $U^i$ are identity matrices, and the capacity constraints corresponding to $e\in E$ with $|I(e)|\ge 2$ form the matrix $[\bar A^1\;\ldots\;\bar A^r]$ where all the $\bar A^i$ are subsets of the rows of identity matrices. Since $\begin{bmatrix}A^i\\ \bar{A}^i
\end{bmatrix}$ is the constraint matrix for flow $i$, this matrix is TU.

Applying the affine TU decomposition from part~(i), we only need integral variables for the arcs that are allowed to be shared by multiple flows. These are $\sum_{e\in E:|I(e)|\ge 2}|I(e)|\le r\cdot\left|\{e\in E\mid |I(e)|\ge 2\}\right|$ integral variables.
\end{enumerate}
\end{examplex}

Next we relate the class of almost totally unimodular matrices to the matrices with affine TU-dimension 1 (called nearly TU matrices by~\cite{Gijswijt2005}).

A square matrix $A$ is called \textit{almost totally unimodular} if $A$ is not TU but every proper submatrix of $A$ is TU. Almost totally unimodular matrices were introduced by Padberg~\cite{Padberg1988} and are studied as a building block for $k$-balanced matrices~\cite[and references therein]{Conforti2006}.

\begin{examplex}
\label{ex:almost-TU}
Assume $A=[a\;\;\bar{A}]\in\{0,\pm 1\}^{n\times n}$ is almost totally unimodular. Then one can write $A=[\0_n\;\bar{A}]+a[1,0,\ldots,0]$. This is a $1$-row affine TU decomposition, since every $(n-1)\times (n-1)$ submatrix of $\bar{A}$ is TU. Thus, every almost totally unimodular matrix has affine TU-dimension 1.  
\end{examplex}

This example shows that integer linear programs described with an almost totally unimodular constraint matrix can be solved in polynomial time by solving a mixed integer linear program with only one integer variable.   The authors of~\cite{veselov-gribanov-2015} study a more general class of matrices called \emph{almost unimodular matrices}.  Through studying the lattice width of polyhedra with constraint matrices that have small subdeterminants, they show that integer linear programs described by almost unimodular matrices can also be solved in polynomial time.

We conclude this section with a discussion of the so-called \emph{integer decomposition property}.  A polyhedron $P$ has the integer decomposition property, if for any positive integer $k$ every integral vector in $kP$ is the sum of $k$ integral vectors from $P$.  A matrix $A\in\Z^{m\times n}$ is TU (has affine TU-dimension 0) if and only if $\{x\in\R^n\mid Ax\le b,\,x\ge \0_n\}$ has the integer decomposition property for each $b\in\R^m$~\cite{Baum1978}. It is straightforward from the proof in \cite{Baum1978} that for a TU matrix $A$ and $b\in\Z^m$, also $P=\{x\in\R^n\mid Ax\le b\}=\conv(\{x\in\Z^n\mid Ax\le b\})$ has the integer decomposition property. Gijswijt~\cite{Gijswijt2005} showed that for $A\in\{0,\pm 1\}^{m\times n}$ with affine TU-dimension $1$ and $b\in\Z^m$, the polyhedron $P_{A,b}:=\conv(\{x\in\Z^n\mid Ax\le b\})$ has the integer decomposition property, too.

The integer decomposition property does not hold in general for polytopes $P_{A,b}$ where the matrix $A$ has affine TU-dimension 2 and $b\in\Z^m$. As a counterexample, consider the parity polytope $P=\conv\left(\{(0,0,0)^T,(0,1,1)^T,(1,0,1)^T,(1,1,0)^T \}\right)$. It does not have the integer decomposition property, since $(1,1,1)^T\in2P\cap\Z^3$ cannot be written as the sum of two points in $P\cap\Z^3$ (see, e.g., \cite{Hibi2012}). The polytope $P$ has a natural inequality description as $P=\{x\in\R^3\mid Ax\le b\}$ with
$$
A=\begin{bmatrix}1 &1 &1\\-1&1&1\\1&-1&1\\1&1&-1\end{bmatrix}\text{ and }b=\begin{bmatrix}2\\0\\0\\0\end{bmatrix}.
$$
$A$ has the 2-row affine TU decomposition
$$
\begin{bmatrix}1 &1 &1\\-1&1&1\\1&-1&1\\1&1&-1\end{bmatrix}=
\begin{bmatrix}1 &0 &0\\-1&0&0\\1& 0&0\\1&0& 0\end{bmatrix}+
\begin{bmatrix}1 &1\\1&1\\-1&1\\1&-1\end{bmatrix}
\begin{bmatrix}0 &1 &0\\0&0&1\end{bmatrix}.
$$
Since $P$ does not have the integer decomposition property, there exists no 1-row affine TU decomposition for $A$. Thus, $A$ has affine TU-dimension 2. In particular, not every $P_{A,b}$ arising from a matrix $A$ with affine TU-dimension 2 has the integer decomposition property.


\section{Determining the TU-dimension}
\label{sec:TU-dim_determining}


\label{subsec:Bounds-on-TUdim}
Seymour~\cite{Seymour1980} gave a decomposition procedure to recognize if a matrix is TU. A careful implementation of this decomposition procedure results in an algorithm that decides in polynomial time if a matrix is TU, see Truemper~\cite{Truemper1990}. Thus one can decide in polynomial time if a matrix has affine TU-dimension 0.
Given any specific value $k$, we would like to decide if $A$ has affine TU-dimension $k$. 

For $B \in\R^{m \times n}$, let $\im(B) := \{B y \in \R^m \mid y \in \R^n\}$ be the column space of $B$ and $\ker(B)=\left\{ x\in\R^{n}\mid Bx=\0_m\right\}$ be the kernel of $B$. 
A TU decomposition of a matrix $A\in\Z^{m\times n}$ is strongly related to the existence of $\pm 1$ combinations of column vectors of $A$ generating the vector $\0_m$, i.e. vectors $r\in\ker(A)\cap \{0,\pm 1\}^n\setminus\{\0_n\}$.

We will need the following fact, a proof of which can for example be found in \cite{Karzanov1997}. 

\begin{lemmax}[Totally unimodular matrices have totally unimodular kernels]
\label{lem:tu-kernel}
Let $A\in\left\{ 0,\pm1\right\}^{\left(n-k\right)\times n}$ be TU and with rank $n-k$.
Then one can find a TU matrix $W\in\left\{ 0,\pm1\right\} ^{k\times n}$ such that $\ker(A) = \im(W^T)$.
\end{lemmax}

We now show that it is \textit{NP}-hard to decide if a matrix has affine TU-dimension $n$.

\begin{theoremx}\label{thm:TUdimension_n_NPhard}
It is an \textit{NP}-hard problem to decide if a matrix $A\in\Z^{m\times n}$ has affine TU-dimension equal to $n$, even when we restrict to positive matrix entries and $m=1$.
\end{theoremx}

\begin{proofx}

It was shown in~\cite{Woeginger1992} that the equal-sum-subsets problem is \textit{NP}-complete. Given $b\in\Z^n_+$, this problem is to decide if there exists a nonzero $r\in\{0,\pm 1\}^n$ such that $b^Tr=0$.

We reduce the equal-sum-subsets problem to deciding if the transpose of the vector $2n\cdot b$ admits a $(n-1)$-row affine TU decomposition.
In other words, for any vector $b\in\Z^n_+$, $a^T=2n\cdot b^T$, we show that $b$ is a yes-instance to the equal-sum-subset problem if and only if there exists a decomposition $a^T=\tilde{a}^T+u^TW$ with $u\in\Z^{n-1}$, $W\in\{0,\pm 1\}^{(n-1)\times n}$, $\tilde{a}\in\{0,\pm 1\}^n$, and $W$ together with the row vector $\tilde{a}^T$ is TU.

If there exists an affine TU decomposition for $a^T$ as above with $W$ TU having rank $k\le n-1$, then by Lemma \ref{lem:tu-kernel} there exists a TU matrix $R\in\{0,\pm 1\}^{(n-k)\times n}$ having rank $n-k\ge 1$ such that $$\{x\in\R^n\mid Wx=\0_{n-1}\}=\{R^Tt\mid t\in\R^{n-k}\}.$$ In particular $R$ has only nonzero rows. Consider its first row $r\in\{0,\pm 1\}^{n}$, which, as all other rows of $W$, satisfies  $Wr=\0_{n-1}$. Then $a^Tr=\tilde{a}^Tr+u^TWr=\tilde{a}^Tr$, thus
$$0=a^Tr-\tilde{a}^Tr=2n\cdot b^Tr-\tilde{a}^Tr\;.$$
Since $\tilde{a},r\in\{0,\pm 1\}^n$, $\tilde{a}^Tr\in\{-n,\ldots,n\}$. Furthermore
$2n\cdot b^Tr\in 2n\Z$, therefore $\tilde{a}^Tr=0$ and
$b^Tr=0$ for the nonzero $r\in\{0,\pm 1\}^n$.

Now assume there is a solution to the equal-sum-subsets problem, a nonzero $r\in\{0,\pm 1\}^n$ such that $b^Tr=0$.
By Lemma \ref{lem:tu-kernel} there exists a TU matrix $W\in\{0,\pm 1\}^{(n-1)\times n}$ such that
$$\{x\in\R^n\mid r^Tx=0\}=\{W^Tu\mid u\in\R^{n-1}\}.$$
Thus $b=W^Tu$ has a solution $u\in\R^{n-1}$, and since $W$ is TU we may choose $u$ to be integral.
Then $a^T=2n\cdot b^T$ admits the $(n-1)$-row affine TU decomposition $a^T=\0_n^T+(2n\cdot u)^TW$, where $W$ together with the row $\0_n^T$ is TU.

Therefore the problem to decide if a given matrix has affine TU-dimension $n$ is \textit{NP}-hard.
\end{proofx}

In view of the hardness result, it is an open question if there exists a polynomial time algorithm to decide if $A$ has affine TU-dimension $k$ when $k$ is fixed. Below we provide polynomial time algorithms for other certain special cases of the task.
We first present some auxiliary facts. A proof of the first one can be found in~\cite[Section 19.3/19.4]{schrijver_theory_1986}.
\begin{lemmax}\label{lem:pivoting}
For a nonsingular matrix $E$,  
$\begin{bmatrix}
B & D\\E & C\end{bmatrix}
$ is TU if and only if  $
\begin{bmatrix}
BE^{-1} & D - BE^{-1}C\\-E^{-1} & E^{-1}C
\end{bmatrix}
$ is TU.
\end{lemmax}

\begin{lemmax}
\label{lem:unimod-transform-yields-W-with-unitmatrix}
Suppose $A\in\Z^{m\times n}$ admits an affine TU decomposition with $W\in\{0,\pm 1\}^{k\times n}$.
Then there exists  a matrix $W' \in \{0, \pm 1\}^{k \times n}$ that contains $I_{k}$ as a submatrix such that $A$ also admits an affine TU decomposition with $W'$.

Furthermore, there exists a unimodular matrix $T \in \Z^{k \times k}$ such that $W' = T W$.  The matrices $T$ and $W'$ can be computed in polynomial time.
\end{lemmax}

\begin{proofx}
Without loss of generality, assume that $W$ has rank $k$.
By reordering the columns of $A$ and $W$ with the same column permutation we may assume that $W_1$, the first $k$ columns of $W=[W_1\;\;W_2]$, have column rank $k$.

Since $W_1$ is unimodular, $W_1^{-1}\in\Z^{k\times k}$ is also unimodular. Then we have $W' :=W_1^{-1} W=[I_k\;\; W'_2]$ with $W'_2=W_1^{-1}W_2\in\Z^{k\times (n-k)}$.

Consider an affine TU decomposition $A=\tilde{A}+UW$, where we define $\tilde{A}=[\tilde{A}_1\;\;\tilde{A}_2]$ with $\tilde{A}_1$ being the first $k$ columns of $\tilde{A}$. In particular, the matrix $\begin{bmatrix}
\tilde{A}_1&\tilde{A}_2\\W_1&W_2
\end{bmatrix}$ is TU. Applying Lemma~\ref{lem:pivoting} with $E=W_1$ shows that $\begin{bmatrix}
\tilde{A}_1W_1^{-1}&\tilde{A}_2-\tilde{A}_1W_1^{-1}W_2\\-W_1^{-1}&W_1^{-1}W_2
\end{bmatrix}$ is TU as well, and so is the matrix $\begin{bmatrix}
 0_{m\times k}&\tilde{A}_2-\tilde{A}_1W_1^{-1}W_2\\I_k&W_1^{-1}W_2
\end{bmatrix}$. In particular, also $W'$ is TU. Since we have
$$A=\begin{bmatrix}\tilde{A}_1&\tilde{A}_2\end{bmatrix}+U\begin{bmatrix}W_1&W_2\end{bmatrix}=\begin{bmatrix} 0_{m\times k}&\tilde{A}_2-\tilde{A}_1W_1^{-1}W_2\end{bmatrix}+(\tilde{A}_1+UW_1)\begin{bmatrix}I_k&W_1^{-1}W_2\end{bmatrix},$$
the latter is also an affine TU decomposition of $A$ with the TU matrix $W'$.  This completes the proof of existence for $T := W_1^{-1}$.  The matrices $T$ and $W'$ can be computed in polynomial time using basic linear algebra techniques.
\end{proofx}

\begin{theoremx}
\label{thm:decide-W-with-unitmatrix}
Let $n>k$, $A_1\in\Z^{m\times k}$, $A_2\in\Z^{m\times (n-k)}$, and $W_2\in\{0,\pm 1\}^{k\times (n-k)}$ TU.
Then $[A_1\;\;A_2]$ admits a $k$-row affine TU decomposition with $U\in\Z^{m\times k}$ and $W=[I_k\;\;W_2]$ if and only if $\begin{bmatrix}
A_2-A_1W_2\\
W_2
\end{bmatrix}$ is TU.

Moreover, given $A$ and $W=[I_k\;\;W_2]$, we can find in polynomial time $U$ and $\tilde{A}$ such that $A=\tilde{A}+UW$ is an affine TU decomposition, if it exists.
\end{theoremx}

\begin{proofx}
For the only if part, assume there is a matrix $U$ such that $[A_1\;\;A_2]=\tilde{A}+U[I_k\;\;W_2]$ and $\begin{bmatrix} \tilde{A} \\
\left[I_k\;\; W_2 \right] \end{bmatrix}$ is TU. Applying Lemma~\ref{lem:pivoting} on the TU matrix $\begin{bmatrix} \tilde{A} \\
\left[I_k\;\; W_2 \right] \end{bmatrix}=\begin{bmatrix} A_1-U & A_2-UW_2 \\
I_k & W_2 \end{bmatrix}$ with $E=I_k$ implies that 
$\begin{bmatrix} A_1-U & A_2-A_1W_2 \\
-I_k & W_2 \end{bmatrix}$ is TU as well. In particular, $\begin{bmatrix} A_2-A_1W_2 \\ W_2 \end{bmatrix}$ is TU.

For the if part, define $\tilde{A}=\left[ 0_{m\times k}\;\;A_2-A_1W_2\right]\in\Z^{m\times n}$ and $U=A_{1}\in\Z^{m\times k}$. Then $$\begin{bmatrix} \tilde{A} \\
\left[I_k\;\; W_2 \right] \end{bmatrix}=\begin{bmatrix}
 0_{m\times k}&A_2-A_1W_2\\I_k&W_2 \end{bmatrix}$$
is TU by assumption, and $[A_1\;\;A_2]=\tilde{A}+U[I_k\;\;W_2]$ by construction.  

Thus, in order to find an affine TU decomposition with given $A$ and $W=[I_k\;\;W_2]$ in polynomial time, we apply Truemper's polynomial time algorithm~\cite{Truemper1990} to decide if $\begin{bmatrix}
A_2-A_1W_2\\
W_2
\end{bmatrix}$ is TU. If this is the case, we have an affine TU decomposition $A=\tilde{A}+UW$ with $\tilde{A}=\left[ 0_{m\times k}\;\;A_2-A_1W_2\right]$ and $U=A_{1}$.
Otherwise, there does not exist such a decomposition .
\end{proofx}

\begin{observation}
\label{obs:affineTUdecomp-fixed-n}
Assume that $n$ is fixed and $A\in\Z^{m\times n}$ is given.

For an affine TU decomposition one has $W \in \{0,\pm 1\}^{k \times n}$ and in view of Remark~\ref{rem:W-full-rank} without loss of generality $W$ has rank $k\le n$. Moreover, by Lemma~\ref{lem:unimod-transform-yields-W-with-unitmatrix} one can limit the search to matrices $W$ that contain $I_k$ as a submatrix. 
Enumerating these polynomially many matrices for all $k\le n$, and applying for each of them the efficient algorithm from Theorem~\ref{thm:decide-W-with-unitmatrix}, one can find a matrix $W$ having fewest number of rows such that $A$ admits an affine TU decomposition together with $W$. 
 
Thus, one can compute the affine TU-dimension of $A$ in polynomial time for fixed $n$. 
\end{observation}

Next we show that for fixed $k$, the number of rows of $W$, and fixed $m$, the number of rows of $A$, we can decide efficiently if there exist $U$ and $W$ such that $A$ admits an affine TU decomposition.  

\begin{theoremx}\label{thm:affineTUdim_m_k_fixed}
Suppose $A \in \Z^{m \times n}$.  Then in polynomial time we can decide if $A$ has affine TU-dimension $k$, provided that $m$ and $k$ are fixed.
\end{theoremx}
\begin{proofx}

Let us first reduce the statement to the case that $A$ has distinct columns.

For this, without loss of generality assume the columns of $A$ to be ordered as $A = [A_1\;\;A_2]$,  where $A_1$ has distinct columns and $A_2$ has columns that are repeats of $A_1$. Let $A_1$ have $r$ columns. Consider any affine TU decomposition of $A=\tilde A+UW$ with a TU matrix $[\tilde A\;\;W]$. Let us divide the latter matrix into its first $r$ columns, and the remaining ones 
$$\begin{bmatrix}\tilde A\\W\end{bmatrix}=\begin{bmatrix}\tilde A_1&\tilde A_2\\W_1&W_2\end{bmatrix}.$$
Then we can we can choose columns $\stack{\tilde A'_2}{W'_2}$ as repeats of columns of $\stack{\tilde A_1}{W_1}$ such that with $\tilde A'=[\tilde A_1\,\tilde A'_2]$ and $W'=[W_1\,W'_2]$ we get the (possibly different) affine TU decomposition $A=\tilde A'+UW'$. In particular, $A$ and $A_1$ have the same affine TU-dimension.

From now on assume that $A\in\Z^{m\times n}$ has distinct columns and gives rise to an affine TU decomposition $A=\tilde{A}+UW$. Then the TU matrix $\stack{\tilde{A}}{W}$ with $m+k$ rows also has distinct columns. A TU matrix with $r$ rows and distinct columns has at most $r^2+r+1$ columns (for a reference, see Remark~\ref{rem:master_lower_bound} above). Thus, $n\le (m+k)^2+m+k+1$.

Now with $k$ and $m$, also $n$ is constant. By Observation~\ref{obs:affineTUdecomp-fixed-n}, we thus can compute the affine TU-dimension of $A$ in polynomial time.
\end{proofx} 

We close this section with a conjecture on the complexity of determining the affine TU-dimension of a matrix.

\begin{conjecturex}
Suppose $A \in \Z^{m \times n}$.  Then in polynomial time we can decide if $A$ has affine TU-dimension $k$, provided that $k$ is fixed.
\end{conjecturex}

The above theorems and algorithms can be easily adapted to the consideration of TU decompositions instead of affine TU decompositions. In particular, the proof of Theorem~\ref{thm:TUdimension_n_NPhard}
does also show the \textit{NP}-hardness of recognizing TU-dimension $n$.

\begin{observation}
It is \textit{NP}-hard to decide if a matrix $A\in\Z^{m\times n}$ has TU-dimension equal to $n$, even when we restrict to positive matrix entries and $m=1$.
\end{observation}


\section{Reformulations for knapsack polytopes}
\label{sec:reformulations-knapsack}


In this section we investigate the size of reformulations specific to both $A$ and $b$.   In particular, we focus on bounds for the minimum number of integrality constraints needed to model a general $0$-$1$ knapsack polytope.  
We then provide an example that demonstrates how adding even a single integrality constraint can vastly reduce the number of inequalities needed to describe the integer hull of a knapsack polytope.

\begin{lemmax}[Lower bound]
\label{prop:lower-bound}
Let $m \in \Z_+$ and consider the knapsack polytope 
$$
P = \left\{ \stack{x}{y} \in [0,1]^{2m} \;\middle\vert\;  2(x_1 + y_1) + 2^2(x_2 + y_2) + \dots + 2^m(x_m + y_m)\leq 2^{m+1} - 1\right\}
$$
in dimension $n=2m$. Let $W,W'\in\Z^{k\times m}$ be integral matrices and let $$Q = \conv\left(P \cap \left\{ \stack{x}{y} \;\middle\vert\; [W\;\;W']\stack{x}{y}  \in \Z^k\right\}\right).$$  If $k < m$, then $Q$ is not an integral polyhedron.
\end{lemmax}

\begin{proofx}

Since $k < m$, $\dim(\ker(W)) \geq 1$.  Let $r \in \ker(W) \setminus \{\0_m\}$ such that $a^T r \geq 0$ where $a^T:=(2,2^2,\ldots,2^m)$.  
Let $x^0 \in \{0,1\}^m$  with $x^0_i=\begin{cases} 0 &\mbox{if } r_i \ge 0 \\ 1 & \mbox{if } r_i <0 \end{cases}$. Then there exists a $\hat \lambda > 0$ with $x^0 + \lambda r \in [0,1]^m$ for all $0 < \lambda < \hat \lambda$.
Let $\overline{x^0}$ be the complement of $x^0$, that is, $\overline{x^0} \in \{0,1\}^m$ and $\overline{x^0} + x^0 = \1_m$.  Observe that $\stack{x^0}{\overline{x^0}} \in P\cap \{0,1\}^{2m}$ since 
$$
a^Tx^0 + a^T \overline{x^0} = a^T \1_m = 2^{m+1} - 2 < 2^{m+1} -1.
$$
Let $0 < \tilde \lambda < \min\{ \hat\lambda, \frac{1}{a^T r}\}$  where if $a^T r = 0$ we set $\frac{1}{a^T r } $ to $\infty$. 

Now, suppose that $Q$ is an integral polytope.  This implies that $Q \cap \left\{\stack{x}{y} \;\middle\vert\; y = \overline{x^0}\right\}$ is an integral polytope, and hence, its projection onto the first $m$ variables,  $Q_{\overline{x^0}} := \left\{ x \;\middle\vert\; \stack{x}{\overline{x^0}}  \in Q \right\} \subseteq [0,1]^m$ is also an integral polytope.  
By the choice of $r$ and $\tilde \lambda$, we have $a^Tx^0 \leq a^T (x^0 + \tilde \lambda r)$ and $$a^T(x^0 + \tilde \lambda r)+a^T\overline{x^0}=\sum_{i=1}^m2^i+\tilde \lambda a^Tr<2^m-1,$$ thus $x^0 + \tilde \lambda r \in Q_{\overline{x^0}}$.  Consider the optimizers $G = \argmax\{a^T x \mid x \in Q_{\overline{x^0}}\}$ and $F = \argmax\{ a^T x \mid x \in Q_{\overline{x^0}} \cap \{0,1\}^m\}$.  Since $Q_{\overline{x^0}}$ is integral, $G = \conv(F)$

By choice of $a$, the value $a^T x$ is distinct for all $x \in \{0,1\}^m$ and by construction of $Q_{\overline{x^0}}$, $F = \{ x^0\}$ and $G = \conv(F) = \{ x^0\}$.  But since $a^Tx^0 \leq a^T (x^0 + \tilde \lambda r)$, this implies that $x^0 + \tilde \lambda r \in G$ and hence $x^0 + \tilde \lambda r = x^0$, which is a contradiction with the fact that $\tilde \lambda > 0$ and $r \neq \0_m$.  

Thus, we conclude that $Q$ was not integral.  
 \end{proofx}

We next provide a positive result for modeling knapsack polytopes and show that every $0$-$1$ knapsack polytope can be modeled using at most $n-2$ integrality constraints together with the upper and lower bounds on the variables.      We first prove a theorem about separation of disjoint polytopes.

\begin{theoremx}\label{thm:n-half-dimensional-face}
Let $P,Q \subseteq \R^n$ be polytopes such that $P \cap Q  = \emptyset$ and $\dim(\conv(P\cup Q)) = n$.  Then there exist $h \in \R^n$, $\alpha_1, \alpha_2 \in \R$ with $\alpha_1 < \alpha_2$ such that $P \subseteq \{x\in\R^n \mid h^Tx \leq \alpha_1\}$, $Q \subseteq \{x\in\R^n \mid h^Tx \geq \alpha_2\}$, and $\max\{ | \{ x \in \verts(P) \mid h^Tx = \alpha_1\}| ,|\{x\in  \verts(Q) \mid h^Tx = \alpha_2\}| \} \geq \ceil{\tfrac{n+1}{2}}$.
\end{theoremx}

\begin{proofx}

If either $P$ or $Q$ is empty, then the theorem is immediate. Assuming $Q=\emptyset$, one may take a valid inequality $h^Tx \leq \alpha_1$ for which $\{x\in P\mid h^Tx = \alpha_1\}$ is a facet of $P$, and then set $\alpha_2 > \alpha_1$.  

So suppose $P,Q \neq \emptyset$.  Let $\hat x \in \relint(P), \hat y \in \relint(Q)$.  Let $\bar P = P + t$ and $\bar Q = Q + t$ where $t = \hat x - 2 \hat y$.  Let $\bar x = \hat x + t$ and $\bar y = \hat y + t$.  Then $\bar x \in \relint(\bar P)$, $\bar y \in \relint(\bar Q)$ and $\bar x = 2 \bar y$.  We will prove the result for $\bar P$ and $\bar Q$, which implies the result for $P$ and $Q$.

Since $\bar P \cap \bar Q = \emptyset$, by the (strict) hyperplane separation theorem, there exists a pair $(\bar h, \bar \alpha)$ such that $\bar h^T x < \bar \alpha$ for all $x \in \bar P$ and $\bar h^T y > \bar \alpha$ for all $y \in \bar Q$.  
Define
$$
H = \{ (h,\alpha) \in \R^n \times \R \mid h^T x \leq \bar \alpha, h^Ty \geq \alpha, \alpha \geq \bar \alpha \:\:\: \forall\, x \in \verts(\bar P), y \in \verts(\bar Q)\}.
$$
Since $(\bar h, \bar \alpha)$ was a strict separator and since $\verts(\bar P)$ and $\verts(\bar Q)$ are finite, there exists an $\epsilon > 0$ sufficiently small such that the point $(\bar h, \bar \alpha + \epsilon)$ strictly satisfies all inequalities of $H$.  Hence,  $(\bar h, \bar \alpha + \epsilon) \in \intr(H)$, and therefore $H$ is full-dimensional.

We claim that $H$ is a bounded polyhedron.  Indeed, let $(\bar r, \bar \beta) \in \rec(H)$ where 
$$\rec(H) = \{(r,\beta)\in \R^n \times \R \mid r^T x \leq 0, r^T y \geq \beta, \beta\geq 0 \:\:\:  \forall \,x \in \verts(\bar P), y \in \verts(\bar Q)\}$$
denotes the recession cone of $H$.  By convexity, $\bar r^T \bar x \leq 0$ and $\bar r^T \bar y \geq \bar\beta \geq 0$.  Since $\bar x = 2 \bar y$, this implies that $\bar r^T \bar x = \bar r^T \bar y = 0$ and $\bar \beta = 0$.  
Since $\bar x \in \relint(\bar P)$ and $\bar y \in \relint(\bar Q)$, we have that $\bar r^Tx = 0$ for all $x \in \bar P$ and $\bar r^T y = 0$ for all $y \in \bar Q$.  Since $\dim(\conv(\bar P \cup \bar Q)) = n$, there exist $z^1, z^2 \in \conv(\bar P \cup \bar Q)$ and a $\lambda > 0$ such that $\lambda \bar r  = z^1 - z^2$.  By definition of $\conv(\bar P \cup \bar Q)$, there exist $x^1, x^2 \in \bar P$, $y^1, y^2 \in \bar Q$, $\mu_1, \mu_2 \in [0,1]$ such that $z^i = \mu_i x^i + (1- \mu_i) y^i$  for $i=1,2$.  Combining this, we see that 
$$
\lambda \| \bar r\|_2^2 = \lambda \bar r^T\bar r = \lambda \bar r^T(\mu_1 x^1 + \mu_2 x^2 + (1- \mu_1) y^1 + (1-\mu_2)  y^2) = 0,
$$
where the last equality comes from the fact that $\bar r^Tx = 0$ for all $x \in \bar P$ and $\bar r^T y = 0$ for all $y \in \bar Q$.  Hence $\bar r = \0_n$ and $\rec(H) = \{(\0_n,0)\}$. Therefore $H$ is bounded (see, e.g., ~\cite[Section 8.2 (5) (iii)]{schrijver_theory_1986}).

Finally, since $H$ is a bounded, full dimensional polytope, there exists a vertex $(\hat h, \hat \alpha)$ of $H$ such that the inequality $\alpha \geq \bar \alpha$ is not tight.  
Since this vertex is defined by $n+1$ tight inequalities from $H$ that are not the inequality $\alpha \geq \bar \alpha$, there exists a set of $n+1$ points corresponding to the tight inequalities, each of them being a vertex of either $\bar P$ or $\bar Q$.  By pigeonhole principle, there must be at least $\ceil{\tfrac{n+1}{2}}$ such vertices in either $\bar P$ or $\bar Q$.  This completes the proof.
\end{proofx}

Consider such a face $F$ of $P$ respectively $Q$ spanned by at least $\left\lceil \frac{n+1}{2}\right\rceil$ vertices. We want to find a TU matrix $W\in\{0,\pm 1\}^{k\times n}$ such that $Wx=d$ with the same $d\in\Z^k$ for all $x\in F$. Below we prove that for all $0$-$1$ knapsack polytopes this is possible for $n\geq 4$ with $k=n-2$, where $P=\conv{\left\{ x\in\left\{0,1\right\}^{n}\mid a^{T}x\le b\right\}}$ and $Q=\conv(\{0,1\}^n\setminus P)$.  Thus, for all $0$-$1$ knapsack polytopes, it suffices to introduce $n-2$ integrality constraints to satisfy property~\eqref{eq:Wprop}.    

\begin{theoremx}
\label{thm:upper-bound}
Let $n\ge4$, $a\in\Z^n$, $b\in\Z$ and $P=\left\{ x\in\left[0,1\right]^{n}\mid a^{T}x\le b\right\} $.
There exists $W\in\Z^{(n-2)\times n}$ such
that \eqref{eq:Wprop} holds.
\end{theoremx}
\begin{proofx}
Let $S=P\cap \{0,1\}^n$ and $S^c=\{0,1\}^n\setminus S$. Without loss of generality we may assume $S\neq\emptyset$ and $S^c\neq\emptyset$.

Since the polytopes $\conv(S)$ and $\conv(S^c)$ are disjoint and $\conv(S \cup S^c)=\left[0,1\right]^{n}$, by Theorem~\ref{thm:n-half-dimensional-face} there exists $h\in\R^n$, $\alpha_1,\alpha_2$ with $\alpha_1<\alpha_2$ such that $\conv(S) \subseteq \{x\in\R^n \mid h^Tx \leq \alpha_1\}$, $\conv(S^c) \subseteq \{x\in\R^n \mid h^Tx \geq \alpha_2\}$, and $| \{ x \in S \mid h^Tx = \alpha_1\}|\ge \ceil{\tfrac{n+1}{2}}$  or $|\{x\in  S^c \mid h^Tx = \alpha_2\}| \geq \ceil{\tfrac{n+1}{2}}$. 

Since $n\ge 4$, $\ceil{\tfrac{n+1}{2}}\ge 3$. In case that $| \{ x \in S \mid h^Tx = \alpha_1\}|\ge \ceil{\tfrac{n+1}{2}}$, let $x^1,x^2,x^3\in\{x \in S \mid h^Tx = \alpha_1\}$ be disjoint vertices  of $S$. Otherwise, choose $x^1,x^2,x^3\in\{x \in S^c \mid h^Tx = \alpha_2\}$ to be disjoint vertices of $S^c$. Consider
$$A=\begin{bmatrix}
(x^2-x^1)\\
(x^3-x^1)
\end{bmatrix}\in\{0,\pm 1\}^{2\times n}.$$
Since $x^1,x^2,x^3\in\{0,1\}^n$, the columns of $A$ are from the set  $\{[0,0]^T,[0,\pm 1]^T,[\pm 1,0]^T, \pm[1,1]^T\}$. Hence, $A$ is TU. Furthermore, the affine hull $\aff(x^1,x^2,x^3)$ satisfies
$$\aff(x^1,x^2,x^3)=\{x^1+A^Ty\mid y\in\R^2\}=\{x^1+x\mid x\in\R^n, Wx=\0_{n-2}\},$$
where the rows of $W$ span the kernel of the TU matrix $A$. By Lemma~\ref{lem:tu-kernel} we may choose $W\in\Z^{(n-2)\times n}$ to be TU. Since $W$ it TU, the set $Q_z=\{x\in[0,1]^n\mid Wx=z\}$ is an integral polytope for every $z\in\Z^{n-2}$. By construction of $W$, for every $z\in\Z^{n-2}$ there exists an $\alpha\in\R$ such that $Q_z\subseteq\{x\in\R^n\mid h^Tx=\alpha\}$. Thus either $Q_z\subseteq \conv(S)$ or $Q_z\subseteq \conv(S^c)$. Therefore $P\cap Q_z=\conv(P\cap\Z^n)\cap Q_z$ and
$$\conv(P\cap \Z^n)=\conv\left(\bigcup_{z\in\Z^{n-2}} P\cap Q_z\right)=\conv(\{x\in P\mid Wx\in\Z^{n-2}\}).$$ 
\end{proofx}

It is an interesting open question to determine the minimum number of integrality constraints needed to describe any knapsack polytope.  We now know this number lies between $\tfrac{n}{2}$ and $n-2$.

Allowing the addition of polynomially many linear inequalities to the mixed integer reformulation is natural and likely advantageous in many scenarios. In this vein, the parity polytope in Section~\ref{sec:Introdcution} has been formulated using a single integrality condition together with the bounds on the variables and one additional inequality. We conclude this section with a family of knapsack examples in dimension $n$, where the polyhedral description consists of $\Theta(n^{k})$ linear constraints. Here, $k$ is a parameter controlling the weight distribution of the instances.  Each of these knapsack polytopes can be described by only $\mathcal{O}(n)$ linear constraints together with one integrality constraint.

\begin{examplex}[Knapsack with big and small weights]
\label{exa:knapsack}
Let $n\in\Z$, $k\in\Z$, $k\ge 3$ and $b\in\R_+$. Consider weights $a_i\in\R_+$ where $\frac{b}{k+1}<a_{i}\le\frac{b}{k}$ for all $i \in S =\{1,\dots,s\}$ and $\frac{k-1}{k+1}b<a_{i}\le b$ for all $i \in B = \{s+1,\dots,n\}$. To simplify notation, assume that $\frac{b}{k+1}<a_{1}\le a_{2}\le\ldots\le a_{s}\le\frac{b}{k}$
and $s\ge k$. We consider $P=\left\{ x\in\left[ 0,1\right] ^{n}\mid a^{T}x\le b\right\} $.

One can show that the facet description of the polytope $\conv(P \cap \Z^n)$ is given by 
\begin{align*}
 x&\ge 0\\
\sum_{j\in B}x_{j}&\le1\\
\sum_{j=i}^{s}x_{j}+\sum_{\substack{j\in B:\\
a_{j}>b-a_{i}
}
}x_{j}+\left(k-1\right)\cdot\sum_{j\in B}x_{j}&\le k & \text{ for all $i=1,\dots,s-k+1$}\\
\sum_{j\in R}x_{j}+\sum_{\substack{j\in B:\\
a_{j}>b-a_{i\left(R\right)}
}
}x_{j}+\left(\left|R\right|-1\right)\cdot\sum_{j\in B}x_{j} & \le\left|R\right| & \text{for all $R\subseteq S$, $1\le \left|R\right|\le k-1,$}
\end{align*}
where $i\left(R\right)$ is the index of an item from $R$ having smallest weight $a_i$.

It follows easily that $\conv(P \cap \Z^n)$ is the convex hull of all $x\in\left[0,1\right]^{n}$ that satisfy 
\begin{align*}
\sum_{j=i}^{s}x_{j}+\sum_{\substack{j\in B:\\
a_{j}>b-a_{i}
}
}x_{j}+\left(k-1\right)\cdot\sum_{j\in B}x_{j}&\le k & \text{for all $i=1,\dots,s$}\\
\sum_{j\in B}x_{j}&\in\left\{ 0,1\right\}.
\end{align*}

Notice that the latter mixed integer description in fact is an affine TU decomposition. However, the above polyhedral description for the considered knapsack polytope is specific for the right-hand-side $b$ and the fixed $k$.
\end{examplex}

\section*{Acknowledgement}
We thank Santanu S. Dey for discussing his idea for the lower bound in Example~\ref{rem:master_lower_bound}. We owe thanks to Shmuel Onn who made us aware of a much simplified version of the proof of Theorem~\ref{thm:affineTUdim_m_k_fixed}. 

We also want to express our gratitude to two anonymous reviewers. Their detailed comments and suggestions on an earlier version of the manuscript led to enhancements on the general structure of our paper, as well as greatly improved the paper in many ways.

\bibliographystyle{abbrv}
\bibliography{references}{}

\providecommand\CheckAccent[1]{\accent20 #1}
\begin{thebibliography}{10}

\bibitem{Balas1998}
E.~Balas.
\newblock Disjunctive programming: Properties of the convex hull of feasible
  points.
\newblock {\em Discrete Applied Mathematics}, 89(1-3):3--44, 1998.

\bibitem{Baum1978}
S.~Baum and L.~E. Trotter.
\newblock {\em Optimization and Operations Research: Proceedings of a Workshop
  Held at the University of Bonn, October 2--8, 1977}, chapter Integer Rounding
  and Polyhedral Decomposition for Totally Unimodular Systems, pages 15--23.
\newblock Springer Berlin Heidelberg, 1978.

\bibitem{carrkonjevod2004}
R.~D. Carr and G.~Konjevod.
\newblock Polyhedral combinatorics.
\newblock In H.~Greenberg, editor, {\em Tutorials on Emerging Methodologies and
  Applications in Operations Research}, chapter~2, pages 1--48. Springer, 2004.

\bibitem{Conforti2006}
M.~Conforti, G.~Cornu\'ejols, and K.~Vu\v{s}kovi\'{c}.
\newblock Balanced matrices.
\newblock {\em Discrete Mathematics}, 306(19-20):2411--2437, 2006.

\bibitem{Eggan198671}
L.~Eggan and M.~Plantholt.
\newblock The chromatic index of nearly bipartite multigraphs.
\newblock {\em Journal of Combinatorial Theory, Series B}, 40(1):71--80, 1986.

\bibitem{Gijswijt2005}
D.~Gijswijt.
\newblock Integer decomposition for polyhedra defined by nearly totally
  unimodular matrices.
\newblock {\em SIAM Journal on Discrete Mathematics}, 19(3):798--806, 2005.

\bibitem{Hassin2004}
R.~Hassin and A.~Levin.
\newblock An efficient polynomial time approximation scheme for the constrained
  minimum spanning tree problem using matroid intersection.
\newblock {\em SIAM Journal on Computing}, 33(2):261--268, 2004.

\bibitem{heller1957}
I.~Heller.
\newblock On linear systems with integral valued solutions.
\newblock {\em Pacific J. Math.}, 7(3):1351--1364, 1957.

\bibitem{Hibi2012}
T.~{Hibi}, A.~{Higashitani}, L.~{Katth{\"a}n}, and R.~{Okazaki}.
\newblock Normal cyclic polytopes and cyclic polytopes that are not very ample.
\newblock {\em Journal of the Australian Mathematical Society}, 96:61--77,
  2014.

\bibitem{Jeroslow1975}
R.~Jeroslow.
\newblock On defining sets of vertices of the hypercube by linear inequalities.
\newblock {\em Discrete Mathematics}, 11(2):119--124, 1975.

\bibitem{Kaibel2013}
V.~Kaibel and K.~Pashkovich.
\newblock Constructing extended formulations from reflection relations.
\newblock In M.~J\"unger and G.~Reinelt, editors, {\em Facets of Combinatorial
  Optimization -- Festschrift for Martin Gr\"otschel}, pages 77--100. Springer
  Berlin Heidelberg, 2013.

\bibitem{Karzanov1997}
A.~V. Karzanov and S.~T. McCormick.
\newblock Polynomial methods for separable convex optimization in unimodular
  linear spaces with applications.
\newblock {\em SIAM Journal on Computing}, 26(4):1245--1275, 1997.

\bibitem{Korte2007}
B.~Korte and J.~Vygen.
\newblock {\em Combinatorial Optimization: Theory and Algorithms}.
\newblock Springer, 4th edition, 2007.

\bibitem{lenstra_integer_1983}
H.~W. Lenstra, Jr.
\newblock Integer programming with a fixed number of variables.
\newblock {\em Mathematics of Operations Research}, 8:538--548, 1983.

\bibitem{lodi-personal}
A.~Lodi.
\newblock Personal communication, 2014 and 2015.

\bibitem{martin1987}
R.~K. Martin.
\newblock Generating alternative mixed-integer programming models using
  variable redefinition.
\newblock {\em Operations Research}, 35(6):820--831, 1987.

\bibitem{Oriolo2013}
G.~Oriolo, L.~Sanit\`{a} , and R.~Zenklusen.
\newblock Network design with a discrete set of traffic matrices.
\newblock {\em Operations Research Letters}, 41(4):390--396, 2013.

\bibitem{Padberg1988}
M.~Padberg.
\newblock Total unimodularity and the {Euler}-subgraph problem.
\newblock {\em Operations Research Letters}, 7(4):173--179, 1988.

\bibitem{schrijver_theory_1986}
A.~Schrijver.
\newblock {\em Theory of Linear and Integer Programming}.
\newblock John Wiley and Sons, New York, 1986.

\bibitem{Seymour1980}
P.~D. Seymour.
\newblock Decomposition of regular matroids.
\newblock {\em Journal of Combinatorial Theory, Series B}, 28(3):305--359,
  1980.

\bibitem{Truemper1990}
K.~Truemper.
\newblock A decomposition theory for matroids. {V. T}esting of matrix total
  unimodularity.
\newblock {\em Journal of Combinatorial Theory, Series B}, 49(2):241--281,
  1990.

\bibitem{50-years-vanderbeck-wolsey}
F.~Vanderbeck and L.~A. Wolsey.
\newblock Reformulation and decomposition of integer programs.
\newblock In M.~J\"{u}nger, T.~M. Liebling, D.~Naddef, G.~L. Nemhauser, W.~R.
  Pulleyblank, G.~Reinelt, G.~Rinaldi, and L.~A. Wolsey, editors, {\em 50 Years
  of Integer Programming 1958-2008}, pages 431--502. Springer, 2010.

\bibitem{veselov-gribanov-2015}
S.~I. Veselov and D.~V. Gribanov.
\newblock On integer programming with almost unimodular matrices and the
  flatness theorem for simplices.
\newblock Manuscript, arXiv:1505.03132 [cs.CG], 2015.

\bibitem{Woeginger1992}
G.~J. Woeginger and Z.~Yu.
\newblock On the equal-subset-sum problem.
\newblock {\em Information Processing Letters}, 42(6):299--302, 1992.

\end{thebibliography}

\end{document}